\documentclass[12pt,a4paper]{article}
\usepackage{amssymb}
\usepackage[utf8]{inputenc}
\usepackage[pdftex]{hyperref}
\usepackage[english]{babel}
\usepackage{amssymb,latexsym,amsmath,amscd,amsthm}
\usepackage{makeidx}
\usepackage{graphicx}
\usepackage{titlefoot}

\newtheorem{Thm}{Theorem}

\newtheorem{Lem}{Lemma}
\newtheorem{Cor}{Corollary}

\newtheorem{Rem}{Remark}

\title{K-theory for the $C^*$-algebras of the solvable Baumslag-Solitar groups}
\author{Sanaz POOYA and Alain VALETTE}

\begin{document}

\maketitle

\baselineskip=16pt

\unmarkedfntext{Keywords : Solvable Baumslag- Solitar groups, K-theory, K-homology, Baum-Connes conjecture, group $C^*$-algebra. MSC classification: 46L80 }

\begin{abstract} We provide a new computation of the K-theory of the group $C^*$-algebra of the solvable Baumslag-Solitar group $BS(1,n)\;(n\neq 1)$; our computation is based on the Pimsner-Voiculescu 6-terms exact \mbox{sequence}, by viewing $BS(1,n)$ as a semi-direct product $\mathbb{Z}[1/n]\rtimes\mathbb{Z}$. We deduce from it a new proof of the Baum-Connes conjecture with trivial coefficients for $BS(1,n)$.
\end{abstract}

\section{Introduction}

Let $G$ be a countable group. The Baum-Connes conjecture with coefficients for $G$ (denoted by $BC_{coef}$), proposes that, for every 
$C^*$-algebra $A$ on which $G$ acts by automorphisms, the analytical assembly map 
$$ 
\mu_{G,A} : KK_i^G(\underbar {EG}, A) \rightarrow K_i(A\rtimes_ r G) \,\,\,\,\,\, i = 0, 1 
$$
is an isomorphism; where $KK_i^G(\underbar{EG}, A)$ denotes the $G$-equivariant K-homology with $G$-compact supports and coefficients in $A$, of the classifying space $\underbar {EG}$ for $G$-proper actions; and
$K_i(A \rtimes_{r} G)$ denotes the analytical K-theory of the reduced crossed product $A \rtimes_{r} G$.
Although $BC_{coef}$ failed to be true in general, it has been proved for several classes of groups. Among them are one-relator groups, see \cite{Oyono}, \cite{Tu}, \cite{BBV}. Furthermore Higson-Kasparov \cite{HK} established $BC_{coef}$ for the class of amenable groups. For $A = \mathbb C$, we come up with the original Baum-Connes conjecture \cite{BC} that was extended by Baum, Connes and Higson  \cite{BCH} to the above stronger formulation.

The solvable Baumslag-solitar groups $BS(1, n) = \langle a, b | a b a^{-1} = b^{n}\rangle$ for $n \in\mathbb{Z}\backslash\{0\}$, 
are both one-relator and amenable, so are located in the \mbox{intersection} of the two above-mentioned classes. Although $BC_{coef}$ is known for them, it seems interesting to provide in this case an explicit description, with explicit generators, of both sides of the conjecture without coefficients. 

So in this paper we give, for the solvable Baumslag-Solitar groups, a direct proof that $\mu_{G,\mathbb{C}}$ is an isomorphism, by identifying the generators of both sides. In Theorem \ref{main} we show that 
 $$K_0(C^*(BS(1,n)))=\mathbb{Z}.[1]$$ and
 $$K_1(C^*(BS(1,n)))=\mathbb{Z}\oplus \mathbb{Z}/|n-1|.\mathbb{Z}, \,\,\, n\neq 1$$
with generators $[a]$ (of infinite order) and $[b]$ (of order $|n-1|$).

To prove this result we view $BS(1,n)$ as a semi-direct by $\mathbb{Z}$, hence $C^*(BS(1,n))$ as a crossed product by $\mathbb{Z}$, and we compute the analytical \mbox{K-groups} of $BS(1, n)$ thanks to the Pimsner-Voiculescu 6-terms exact \mbox{sequence} \cite{PV} \footnote{For $n=-1$, i.e. the Klein bottle group, a computation based on the Pimsner-Voiculescu sequence appears in Proposition 2.1 of \cite{Sudo}, apparently not aware of previous results on the subject.}.

Finally to reach the Baum-Connes conjecture with trivial coefficients for $BS(1,n)$ we appeal to two useful facts: on the one hand for $G$ a torsion-free group we have $K_i^G(\underbar{EG})=K_i(BG)$, the K-homology with compact supports of a classifying space $BG$ for $G$; on the other hand for $G$ one-relator torsion-free, there is a simple 2-dimensional model for $BG$, namely the presentation complex of $G$, see \cite{Lyn}.




\section{The $C^*$-algebra of $BS(1,n)$}

For $n\neq 1$, there is a faithful homomorphism from $BS(1,n)$ to the affine group of the real line, given by:
$$BS(1,n)\rightarrow Aff_1(\mathbb{R}):\left\{\begin{array}{ccll}a & \mapsto & (x\mapsto nx) & \mbox{(dilation by $n$)} \\b & \mapsto & (x\mapsto x+1) & \mbox{(translation by $+1$)}\end{array}\right. .$$
It realizes an isomorphism 
$$BS(1,n)\simeq \mathbb{Z}[1/n]\rtimes_\alpha\mathbb{Z},$$
where $\mathbb{Z}[1/n]=\{\frac{m}{n^\ell}\in\mathbb{Q}:m, n\in\mathbb{Z},\ell\in\mathbb{N}\}$, viewed as an additive group; and $\alpha$ is multiplication by $n$.

It is well-known that, if a discrete group $G$ decomposes as a semi-direct product $G=H\rtimes_\alpha\mathbb{Z}$, with $H$ a normal abelian subgroup, then 
$$C^*(G)=C^*(H)\rtimes_\alpha\mathbb{Z}=C(\hat{H})\rtimes_{\hat{\alpha}}\mathbb{Z},$$
where $\hat{H}$ denotes the Pontryagin dual of $H$ (so $\hat{H}$ is a compact abelian group), and $\hat{\alpha}$ is the dual automorphism.

In the case of $BS(1,n)$, we have $H=\mathbb{Z}[1/n]$, viewed as the inductive limit of 
$$\mathbb{Z}\stackrel{i_0}{\longrightarrow}\mathbb{Z}\stackrel{i_1}{\longrightarrow}\mathbb{Z}\stackrel{i_2}{\longrightarrow}...,$$
where $i_k:\mathbb{Z}\rightarrow\mathbb{Z}$ (for $k\geq 0$) is multiplication by $n$. So $\widehat{i_k}:\mathbb{T}\rightarrow\mathbb{T}$ is raising to the power $n$, and $\hat{H}$ is the projective limit of 
$$...\stackrel{\widehat{i_2}}{\longrightarrow}\mathbb{T}\stackrel{\widehat{i_1}}{\longrightarrow}\mathbb{T}\stackrel{\widehat{i_0}}{\longrightarrow}\mathbb{T},$$
which we identify with the solenoid\footnote{Strictly speaking, it is a solenoid only for $|n|>1$, while it is $\mathbb{T}$ for $|n|=1$.}
$$X_n=\{z=(z_k)_{k\geq 0}\in\mathbb{T}^\mathbb{N}: z_{k+1}^n=z_{k},\forall k\geq 0\}.$$
The duality between $X_n$ and $\mathbb{Z}[1/n]$ is given by $(z,m)=z_\ell^m$, where $m$ belongs to the $\ell$-th copy of $\mathbb{Z}$; this is well defined as $(z,i_\ell(m))=z_{\ell+1}^{n.m}=z_\ell^m=(z,m)$. For $\frac{m}{n^\ell}\in\mathbb{Z}[1/n]$, this corresponds to $(z,\frac{m}{n^\ell})=z_\ell^m$ for $z=(z_k)_{k\geq 0}\in X_n$.

The automorphism $\alpha$ is given by $\alpha(m)=i_\ell(m)$, where $m$ lies in the $\ell$-th copy of $\mathbb{Z}$. So $\hat{\alpha}$ is the automorphism of $X_n$ given by the backwards shift: $(\hat{\alpha}(z))_k=z_{k+1}$ for $k\geq 0$. 

So $C^*(BS(1, n)) = C(X_n) \rtimes_{\hat \alpha} \mathbb Z $. This crossed product can be viewed as the universal $C^*$-algebra generated by two unitaries $u$ and $v$ satisfying the relation $u v u^{-1} = v^n$, where $u$ is the unitary of $C^*(\mathbb Z)$ corresponding to the generator $+1$ of $\mathbb Z$ acting on $C(X_n)$, while $v\in C(X_n)$ is given by the function $z\mapsto z_0$ on $X_n$. This crossed product description of $C^*(BS(1,n))$ appears already in \cite{BJ,IMSS}.






\begin{Lem}\label{K*(Xn)} $K_0(C(X_n))=\mathbb{Z}.[1]$ (the infinite cyclic group generated by the class of $1\in C(X_n)$) and $K_1(C(X_n))\simeq\mathbb{Z}[1/n]$.
\end{Lem}

{\bf Proof:} We have $C(X_n)=C^*(\mathbb{Z}[1/n])=\varinjlim (C^*(\mathbb{Z}),i_k)$ (where we also denote by $i_k$ the $*$-homomorphism $C^*(\mathbb{Z})\rightarrow C^*(\mathbb{Z})$ associated with the group homomorphism $i_k$). Since K-theory commutes with inductive limits, we get $K_i(C(X_n))=\varinjlim (K_i(C^*(\mathbb{Z})),(i_k)_*)$ ($i=0,1$). Since $K_0(C^*(\mathbb{Z}))=\mathbb{Z}.[1]$ and $i_k$ is a unital $*$-homomorphism, we have \mbox{$K_0(C(X_n))=\varinjlim (\mathbb{Z}.[1],Id)=\mathbb{Z}.[1]$}. On the other hand, let $v$ be the unitary of $C^*(\mathbb{Z})$ corresponding to the generator $+1$ of $\mathbb{Z}$ (so that $K_1(C^*(\mathbb{Z}))=\mathbb{Z}.[v]$). Then $i_k(v)=v^n$, i.e. $(i_k)_*[v]=n[v]$, and the inductive system $(K_1(C^*(\mathbb{Z})),(i_k)_*)$ is isomorphic to the original system $(\mathbb{Z},i_k)$, so they have the same limit $\mathbb{Z}[1/n]$.
\hfill$\square$

\section{K-theory for $C^*(BS(1,n))$}

Let $A$ be a unital $C^*$-algebra and $\alpha \in Aut (A)$. We can define the crossed product $A\rtimes_{\alpha} \mathbb{Z}$ associated with the action $\alpha$ of $\mathbb Z$ on $A$. Let $u \in A\rtimes_{\alpha} {\mathbb Z}$ be the unitary which implements this action in the construction of crossed product. Abstractly the 
crossed product $A\rtimes_{\alpha} \mathbb{Z}$ is generated by \mbox{$\{A, u : u a u^* = \alpha(a),  a \in A\}$}.  The Pimsner-Voiculescu 6-term exact sequence \cite{PV} gives us a tool to calculate the $K$-theory of $A\rtimes_{\alpha} \mathbb{Z}$ via the following cyclic diagram with 6-terms:

\begin{equation*}
 \begin{array}{ccccc}
	K_0(A) 	& \stackrel{Id-\alpha_*}{\longrightarrow} & K_0(A) & \stackrel{\iota_*}{\longrightarrow} & K_0(A\rtimes_{\alpha} {\mathbb Z}) \\ \partial_1\uparrow &   &   &   & \downarrow \partial_0\\
	K_1(A\rtimes_{\alpha} {\mathbb Z}) & \stackrel{\iota_*}{\longleftarrow} & K_1(A) & \stackrel{Id-\alpha_*}{\longleftarrow} & K_1(A)
 \end{array}
\end{equation*}

Here $\iota: A\rightarrow A\rtimes_{\alpha}\mathbb{Z}$ denotes inclusion. We will need some \mbox{understanding} of the connecting map $\partial_1$; namely, we observe in the next lemma that $\partial _1 ([u]) = -[1]$. This will help us in later computations.

\begin{Lem}\label{K_1 generator}
The connecting map $\partial_1 \colon K_1(A\rtimes_{\alpha} \mathbb Z) \rightarrow K_0(A)$ maps $[u]$ to $-[1]$.
\end{Lem}

{\bf Proof:} Let $C^*(S)$ be the $C^*$-algebra generated by a non-unitary isometry $S$ and let $P = I - S^*S$. Now $\mathcal T_{\mathrm{A}, {\alpha}}$, the Toeplitz algebra for $A$ and $\alpha$, is the $C^*$-subalgebra of $(A\rtimes_{\alpha} \mathbb Z)\otimes C^*(S)$ generated by $u \otimes S$ and $A\otimes I$. Let $\mathcal K$ be the $C^*$-algebra of compact operators on a separable Hilbert space, with the corresponding system of matrix units $( e_{ij})_{i,j \geq 0}$. Consider the Toeplitz extension associated with $A\rtimes_{\alpha} \mathbb Z $ as in \cite{PV}:
\begin{equation*}
    0
    \rightarrow
    A \otimes \ \mathcal{K}
    \stackrel{\varphi}{\rightarrow}
    \mathcal T_{\mathrm{A}, {\alpha}} \stackrel{\psi}{\rightarrow}
    \mathrm A \rtimes \mathbb Z
    \rightarrow
     0,
\end{equation*}
	with 
	$\varphi (a \otimes e_{ij}) = u^i a u^{*^{j}} \otimes S^i P S^{* ^j}$ and $$\psi (u \otimes S) = u, \hspace{15pt} \psi (a \otimes I) = a$$ for any $a \in A $ and $i, j \in \mathbb N$.

	The map $\partial_1 \colon K_1(A\rtimes_{\alpha} \mathbb Z) \rightarrow K_0(A\otimes\mathcal{K})$ is then the boundary map \mbox{associated} with the Toeplitz extension, we compute $\partial_1([u])$ following the description given in \cite{Black}, 8.3.1. Consider first
$
   \left (
         \begin{array}{cc}
            u & 0 \\
            0 & u^*
          \end{array}
   \right )
$.
This matrix can be lifted via $\psi$ to a matrix 
$
    M =
      \left(
         \begin{array}{cc}
            u\otimes S & 1\otimes P \\
                     0 & u^*\otimes S^*
          \end{array}
       \right)
$, 
where $M\in \mathrm U_{2}(\mathcal{T}_{A, \alpha}) $.
For $p_1 := (1\otimes I) \oplus 0 \in \mathrm M_2(\mathcal T_{A, \alpha})$ we have
\begin{equation*}
       Mp_1M^* - p_1
       =
      (1 \otimes SS^* - 1 \otimes I) \oplus 0
       =
      (-1 \otimes P ) \oplus 0
      .		
\end{equation*}
The pullback of this element via $\varphi$ is
 $z := (-1 \otimes e_{00}) \oplus 0 \in \mathrm M_2(A \otimes \mathcal K).$ So $\partial_1([u])=-[-z]$.
Via the isomorphism
 $K_0(A \otimes \mathcal K) \cong K_0 (A)$
, the element
  $[-z] = [1 \otimes e_{00}]$
corresponds to $[1]$.
Hence
 $\partial_1 ([u]) = -[1].$ 
\hfill$\square$

\begin{Thm}\label{main} $K_0(C^*(BS(1,n)))=\mathbb{Z}.[1]$. For $n\neq 1$:
$$K_1(C^*(BS(1,n)))=\mathbb{Z}\oplus \mathbb{Z}/|n-1|.\mathbb{Z},$$
 with generators $[a]$ (of infinite order) and $[b]$ (of order $|n-1|$).
 \end{Thm}
 
 {\bf Proof:} We view $C^*(BS(1,n))$ as the crossed product $C^*(BS(1,n))=C^*(\mathbb{Z}[1/n])\rtimes_\alpha\mathbb{Z}$, and apply the Pimsner-Voiculescu 6-terms exact sequence to it. Denoting by the $\iota:C^*(\mathbb{Z}[1/n])\rightarrow C^*(BS(1,n))$ the inclusion, and appealing to lemma \ref{K*(Xn)}, we get:
 
$$ \begin{array}{ccccc}\mathbb{Z}.[1] & \stackrel{Id-\alpha_*}{\longrightarrow} & \mathbb{Z}.[1] & \stackrel{\iota_*}{\longrightarrow} & K_0(C^*(BS(1,n))) \\ \partial_1\uparrow &   &   &   & \downarrow \\K_1(C^*(BS(1,n))) & \stackrel{\iota_*}{\longleftarrow} & \mathbb{Z}[1/n] & \stackrel{Id-\alpha_*}{\longleftarrow} & \mathbb{Z}[1/n]\end{array}$$

Since $\alpha(1)=1$, the upper-left arrow is the zero map. The bottom-right arrow is given by multiplication by $1-n$ on $\mathbb{Z}[1/n]$, so it is injective, hence the right vertical arrow is zero. This shows that $\iota_*:\mathbb{Z}.[1]\rightarrow K_0(C^*(BS(1,n)))$ is an isomorphism. 

Turning to $K_1$, we observe that the relation $[b]=[aba^{-1}]=[b^n]$ implies $(n-1).[b]=0$, i.e. the order of $[b]$ divides $|n-1|$. To prove that this is exactly $|n-1|$, we look at the bottom line of the Pimsner-Voiculescu sequence. Since $\mathbb{Z}[1/n]/Im(Id-\alpha_*)=\mathbb{Z}/|n-1|.\mathbb{Z}$, we get a short exact sequence:
$$0\rightarrow\mathbb{Z}/|n-1|.\mathbb{Z}\rightarrow K_1(C^*(BS(1,n)))\stackrel{\partial_1}{\rightarrow} \mathbb{Z}.[1]\rightarrow 0.$$
which splits to give $K_1(C^*(BS(1,n)))=\mathbb{Z}\oplus \mathbb{Z}/|n-1|.\mathbb{Z}$, with $[b]$ a generator of order $|n-1|$. Since $\partial_1 ([a])=-[1]$ by lemma \ref{K_1 generator}, we see that $[a]$ is a generator of infinite order.
\hfill$\square$

\begin{Cor}\label{BC} Set $G_n=:BS(1,n)$. For $n\neq 1$, the Baum-Connes conjecture without coefficients holds for $G_n$, i.e. the Baum-Connes assembly map \mbox{$\mu_{G_n,\mathbb{C}}:K_i(BG_n)\rightarrow K_i(C^*(G_n))\;(i=0,1)$} is an isomorphism.
\end{Cor}

{\bf Proof:} We appeal to a result of Lyndon \cite{Lyn}: for a torsion-free one-relator group $G=<S|r>$ on $m$ generators, the presentation complex (consisting of one vertex, $m$ edges and one 2-cell) is a 2-dimensional model for the classifying space $BG$. By lemma 4 in \cite{BBV}:
$$K_0(BG)=H_0(BG,\mathbb{Z})\oplus H_2(BG,\mathbb{Z})\;and\;K_1(BG)=H_1(BG,\mathbb{Z}).$$
Moreover $H_2(BG,\mathbb{Z})=0$ when $r$ is not in the commutator subgroup of the free group on $S$. This applies to $G_n$, as we assume $n\neq 1$. 

Then $K_0(BG_n)=H_0(BG_n,\mathbb{Z})=\mathbb{Z}$, generated by the inclusion of a base point. By Example 2.11 on p.97 of \cite{MV}, the image of this element under $\mu_{G_n,\mathbb{C}}$ is $[1]$, the class of $1$ in $K_0(C^*(G_n))$. The result for $K_0$ then follows from Theorem \ref{main}.

Now, for any group $G$, identify $H_1(BG,\mathbb{Z})$ with the abelianized group $G^{ab}$. There is a map $\kappa_{G}:G^{ab}\rightarrow K_1(C_r^*(G))$ obtained by mapping a group element $g\in G$ first to the corresponding unitary in $C_r^*(G)$, then to the class $[g]$ of this unitary in $K_1(C_r^*(G))$. We get this way a homomorphism $G\rightarrow K_1(C^*_r(G))$, which descends to $\kappa_G:G^{ab}\rightarrow K_1(C^*_r(G))$ as the latter group is abelian. By Theorem 1.4 on p.86 of \cite{MV}, for $G$ torsion-free, the map $\kappa_G$ coincides with $\mu_{G,\mathbb{C}}$ on the lowest-dimensional part of $K_1(BG)$. Here, $\mu_{G_n,\mathbb{C}}:K_1(BG_n)\rightarrow K_1(C^*(G_n))$ coincides with $\kappa_{G_n}:G_n^{ab}=\mathbb{Z}\oplus\mathbb{Z}/|n-1|.\mathbb{Z}\rightarrow K_1(C^*(G_n))$, which is an isomorphism by Theorem \ref{main}.
\hfill$\square$
\begin{Rem}
	Let $\tau : C^*(BS(1, n)) \rightarrow \mathbb C$ be the canonical trace on 
	$ C^*(BS(1, n))$. This induces the homomorphism  
	$\tau _* : K_0( C^*(BS(1, n))) \rightarrow \mathbb R$ 
	at the K-theory level. Since $\tau$ is unital, by Theorem \ref{main} we have $\tau_*(K_0( C^*(BS(1, n)))) = \mathbb Z$. 
\end{Rem}

\vspace{1cm}
Authors addresses:

Institut de Math\'ematiques - Universit\'e de Neuch\^atel

Unimail - Rue Emile Argand 11

CH-2000 Neuch\^atel - SWITZERLAND

sanaz.pooya@unine.ch

alain.valette@unine.ch


\begin{thebibliography}{X}
 
 \bibitem[BC]{BC} Paul {\sc Baum}, Alain {\sc Connes}, 
 \newblock {\em Geometric K-theory for Lie groups and foliations}, 
 \newblock Enseign. Math. 46 (2000), 3-42

\bibitem[BCH]{BCH} P. {\sc Baum}, A. {Connes}, N. {\sc Higson}
\newblock{\em Classifying space for proper actions and K-theory of group $C^*$-algebras}
\newblock Contemp. Math 164 (1994), 241-292.

\bibitem[BBV]{BBV} C. {\sc B\'eguin}, H. {\sc Bettaieb}, A. {\sc Valette}
\newblock {\em K-theory for the $C^*$-algebras of one-relator groups}
\newblock K-theory 16 (1999), 277-298.

\bibitem[Bla]{Black} B. {\sc Blackadar}
\newblock {\em K-Theory for Operator Algebras},
\newblock Mathematical Sciences Research Institute Publications 5, Springer-Verlag, 1986.

\bibitem[BJ]{BJ} B. {\sc Brenken} and P. E. T. {\sc J\o rgensen}, 
\newblock {\em A family of dilation crossed product algebras}, 
\newblock J. Operator Theory 25 (1991), 299-308.

\bibitem[HK]{HK} N. {\sc Higson}, G. {\sc Kasparov},
\newblock{\em Operator K-theory for groups which act properly and isometrically on Hilbert space}
\newblock Electron. Res. Announc. Amer. Math. Soc. 3 (1997), 131-142.

\bibitem[IMSS]{IMSS} B. {\sc Iochum}, T. {\sc Masson}, T. {\sc Sch\"ucker}, A. {\sc Sitarz}, 
\newblock {\em $\kappa$-Deformation and Spectral Triples}
\newblock Acta Phys. Polon. Supp. 4 (2011), 305-324.

\bibitem[Lyn]{Lyn} R.C. {\sc Lyndon},
\newblock {\em Cohomology theory of groups with a single defining relation}
\newblock Ann. of Math. 52 (1950), 650-665.

\bibitem[MV]{MV} G. {\sc Mislin} and A. {\sc Valette}
\newblock {\em Proper group actions and the Baum-Connes conjecture},
\newblock Advanced Courses in Mathematics - CRM Barcelona, Birkh\"auser, 2003.

\bibitem[Oyo]{Oyono} H. {\sc Oyono-Oyono}
\newblock {\em La conjecture de Baum- Connes pour les groupes agissant sur les arbes}
\newblock C. R. Acad. Sci. Paris, S\'er. I 326 (1998), 799-804.

\bibitem[Tu] {Tu} J.-L. {\sc Tu}
\newblock {\em The Baum-Connes conjecture and discrete group actions on trees}
\newblock K-Theory 17 (1999), 303–318.

\bibitem[PV]{PV} M. {\sc Pimsner}, D. {\sc Voiculescu},
\newblock{\em Exact sequences for K-groups and Ext-groups of certain cross product $C^*$- algebras} 
\newblock Operator theory, 4 (1980), 93-118.

\bibitem[Sud]{Sudo} T. {\sc Sudo}
\newblock {\em K-theory for the group $C^*$-algebras of certain solvable discrete groups}
\newblock Hokkaido Math. J. 43 (2014), 209-260.

\end{thebibliography}
\end{document}